\documentclass[12pt,twoside,openany]{paper}
\usepackage{enumerate,amsmath,graphics,amssymb,graphicx,amscd,xypic,amscd,amsbsy,multirow,float,booktabs,verbatim}
\newcommand{\QED}{\hspace*{\fill}\rule{2.5mm}{2.5mm}}
\newtheorem{theorem}{Theorem}[section]
\newtheorem{definition}{Definition}[section]

\newenvironment{proof}{\noindent{\bf Proof\ }}{\QED\\}
\newcommand{\R}{\mathbb{R}}

\newcommand{\N}{\mathbb{N}}

\newtheorem{lemma}{Lemma}[section]

\newtheorem{corollary}{Corollary}[section]

\begin{document}
\begin{center}
\vspace{0.5cm} {\large \bf ``Divergence of sample quantiles"}\\
\vspace{1cm} Reza Hosseini, Simon Fraser University\\
Statistics and actuarial sciences, 888 University Road,\\
Burnaby, BC, Canada, V65 1S6\\
 reza1317@gmail.com
\end{center}

\begin{abstract}

We show that the left (right) sample quantile tends to the left
(right) distribution quantile at $p \in [0,1]$, if the left and
right quantiles are identical at $p$. We show that the sample
quantiles diverge almost surely otherwise. The latter can be
considered as a generalization of the well-known result that the sum
of a random sample of a fair coin with 1 denoting heads and -1
denoting tails is 0 infinitely often. In the case that the sample
quantiles do not converge we show that the limsup is the right
quantile and the liminf is the left quantile.

\end{abstract}


\section{Introduction}

The traditional definition of quantiles for a random variable $X$
with distribution function $F$,
\[lq_X(p)=\inf \{x|F(x) \geq p\},\]
appears in classic works as \cite{parzen-1979}. We call this the
``left quantile function''. In some books (e.g. \cite{rychlik}) the
quantile is defined as
\[rq_X(p)=\inf \{x|F(x) > p\}=\sup \{x| F(x) \leq p\},\]
this is what we call the ``right quantile function''. Also in
robustness literature people talk about the upper and lower medians
which are a very specific case of these definitions. Hosseini in
\cite{reza-phd} considers both definitions, explore their relation
and show that considering both has several advantages.

Section \ref{section-quantile-def-limit-theory} studies the limit
properties of left and right quantile functions. In Theorem
\ref{limit-lq-rq}, we show that if left and right quantiles are
equal, i.e. $lq_F(p)=rq_F(p)$, then both sample versions
$lq_{F_n},rq_{F_n}$ are convergent to the common distribution value.
We found an equivalent statement in Serfling \cite{serfling} with a
rather similar proof. The condition for convergence there is said to
be $lq_F(p)$ being the unique solution of $F(x-)<p\leq F(x)$ which
can be shown to be equivalent to $lq_F(p)=rq_F(p)$. Note how
considering both left and right quantiles has resulted in a cleaner,
more comprehensible condition for the limits. In a problem Serfling
asks to show with an example that this condition cannot be dropped.
We show much more by proving that if $lq_F(p) \neq rq_F(p)$ then
both $rq_{F_n}(p)$ and $rq_{F_n}(p)$ diverge almost surely. The
almost sure divergence result can be viewed as an extension to a
well-known result in probability theory which says that if
$X_1,X_2,\cdots$ an i.i.d sequence from a fair coin with -1 denoting
tail and 1 denoting head and $Z_n=\sum_{i=1}^n X_i$ then $P(Z_n=0 \;
i.o.)=1.$ The proof in \cite{breiman} uses the Borel-Cantelli Lemma
to get around the problem of dependence of $Z_n$. This is equivalent
to saying for the fair coin both $lq_{F_n}(1/2)$ and $rq_{F_n}(1/2)$
diverge almost surely. For the general case, we use the
Borel-Cantelli Lemma again. But we also need a lemma (Lemma
\ref{lemma-sum-deviations}) which uses the Berry-Esseen Theorem in
its proof to show the deviations of the sum of the random variables
can become arbitrarily large, a result that is easy to show as done
in \cite{breiman} for the simple fair coin example. Finally, we show
that even though in the case that $lq_F(p) \neq rq_F(p)$,
$lq_{F_n},rq_{F_n}$ are divergent; for large $n$s they will fall in
\[(lq_F(p)-\epsilon,lq_F(p)] \cup [rq_F(p),rq_F(p)+\epsilon).\]
In fact we show that
\[\liminf_{n \rightarrow \infty} lq_{F_n}(p)=\liminf_{n \rightarrow \infty} rq_{F_n}(p)=lq_F(p)\]
and
\[\limsup_{n \rightarrow \infty} lq_{F_n}(p)=\limsup_{n \rightarrow \infty} rq_{F_n}(p)=rq_F(p).\]
The proof is done by constructing a new random variable $Y$ from the
original random variable $X$ with distribution function $F_X$ by
shifting back all the values greater than $rq_X(p)$ to $lq_X(p)$.
This makes $lq_Y(p)=rq_Y(p)$ in the new random variable. Then we
apply the convergence result to $Y$.

\section{Limit theory}

\label{section-quantile-def-limit-theory}

To prove limit results, we need some limit theorems from probability
theory that we include here for completeness and without proof.
Their proofs can be found in standard probability textbooks and
appropriate references are given below. If we are dealing with two
samples, $X_1,\cdots,X_n$ and $Y_1,\cdots,Y_n$, to avoid confusion
we use the notation $F_{n,X}$ and $F_{n,Y}$ to denote their
empirical distribution functions respectively.

\begin{definition}
Suppose $X_1,X_2,\cdots,$ is a discrete--time stochastic process.
Let $\mathcal{F}(X)$ be the $\sigma$-algebra generated by the
process and $\mathcal{F}(X_n,X_{n+1},\cdots)$  the $\sigma$-algebra
generated by $X_{n},X_{n+1},\cdots$. Any $E \in \mathcal{F}(X)$ is
called a tail event if $E \in \mathcal{F}(X_n,X_{n+1},\cdots)$ for
any $n \in \N$.
\end{definition}

\begin{definition}

Let $\{A_n\}_{n \in \N}$ be any collection of sets. Then
$\{A_n\;i.o.\}$, read as $A_n$ happens infinitely often is defined
by:

\[\{A_n\;i.o.\}=\cap_{i \in \N}\cup_{j=i}^{\infty} A_j.\]

\end{definition}

\begin{theorem}(Kolmogorov 0--1 law):\\
$E$ being a tail event implies that $P(E)$ is either 0 or 1.
\end{theorem}

\begin{proof}
See \cite{breiman}.
\end{proof}

\begin{theorem}(Glivenko-Cantelli Theorem):\\
 Suppose,
$X_1,X_2,\cdots,$ $i.i.d$, has the sample distribution function
$F_n$. Then
\[\lim_{n \rightarrow \infty} \sup_{x \in \R}|F_n(x)-F(x)| \rightarrow 0,\;\;\;a.s..\]

\end{theorem}
\begin{proof}
See \cite{billin-p}.
\end{proof}

\begin{comment}
Here, we extend the Glivenko-Cantelli Theorem to $F^o,G^o$ and
$G^c$.

\begin{lemma}
Suppose $X$ is a random variable and consider the associated
distribution functions $F_X^o,G_X^o$ and $G_X^c$ with corresponding
sample distribution functions $F_{X,n}^o,G_{X,n}^o$ and $G_{X,n}^c$.
Then

\[\sup_{x \in \R} |G_{X,n}^o-G_X^o| \rightarrow 0,\;a.s.,\]
\[\sup_{x \in \R} |F_{X,n}^o-F_X^o| \rightarrow 0,\;a.s.,\]
and
\[\sup_{x \in \R} |G_{X,n}^c-G_X^c| \rightarrow 0,\;a.s..\]

\label{lemma-gil-cant-dist-fcns}
\end{lemma}

\begin{proof}
 Note that
\[F_{X}^c+G_{X}^o=1 \Rightarrow G_{X}^o=1-F_{X}^c,\]
and
\[F_{X,n}^c+G_{X,n}^o=1 \Rightarrow G_{X,n}^o=1-F_{X,n}^c.\]
Since Glivenko-Cantelli Theorem holds for $F_X^c$ it also holds for
$G_X^o$.

To show the result for $F_X^o$, note that $F_X^o(x)=G^o_{-X}(-x)$
and $F_{X,n}^o(x)=G^o_{-X,n}(-x)$. Also to show the result for
$G_X^c$ note that $G_X^c=1-F_X^o$ and $G_{X,n}^c=1-F_{X,n}^o$.

\end{proof}

\end{comment}

\begin{theorem}(Borel-Cantelli lemma):\\
Suppose $(\Omega,\mathcal{F},P)$ is a probability space. Then
\begin{enumerate}

\item $A_n \in \mathcal{F}$ and $\sum_{1}^{\infty} P(A_n) < \infty
\Rightarrow P(A_n\; i.o)=0$.

\item $A_n \in \mathcal{F}$ independent events with
$\sum_{1}^{\infty} P(A_n) =\infty \Rightarrow P(A_n\; i.o)=1$, where
$i.o.$ stands for infinitely often.

\end{enumerate}

\end{theorem}

\begin{proof}
See \cite{breiman}.
\end{proof}

\begin{theorem}(Berry-Esseen bound):
Let $X_1,X_2,\cdots,$ be $i.i.d$ with $E(X_i)=0<\infty$,
$E(X_i^2)=\sigma$ and $E(|X_i|^3)=\rho$. If $G_n$ is the
distribution of \[{X_1+\cdots+X_n}/{\sigma \sqrt{n}}\] and $\Phi(x)$
is the distribution function of a standard normal random variables
then
\[|G_n(x)-\Phi(x)|\leq 3\rho/\sigma^3\sqrt{n}.\]
\end{theorem}

\begin{corollary}
Let $X_1,X_2,\cdots,$ be $i.i.d$ with $E(X_i)=\mu<\infty$,
$E(|X_i-\mu|^2)=\sigma$ and $E(|X_i-\mu|^3)=\rho$. If $G_n$ is the
distribution of $(X_1+\cdots+X_n-n\mu)/{\sigma
\sqrt{n}}=\sqrt{n}(\frac{\bar{X}_n-\mu}{\sigma})$ and $\Phi(x)$ is
the distribution function of a standard normal random variable then
\[|G_n(x)-\Phi(x)|\leq 3\rho/\sigma^3\sqrt{n}.\]
\end{corollary}
\begin{proof}
This corollary is obtained by applying the theorem to $Y_i=X_i-\mu$.
\end{proof}
Now let $A_n=(X_1+\cdots+X_n-n\mu)/{\sigma \sqrt{n}}$. Then

\[|P(A_n>x)-(1-\Phi(x))|=|P(A_n \leq x)-\Phi(x)|=|G_n(x)-\Phi(x)|<3\rho/\sigma^3\sqrt{n}.\]
Also
\[|P(x<A_n\leq y)-(\Phi(y)-\Phi(x)))| \leq |G_n(y)-\Phi(y)|+|G_n(x)-\Phi(x)| \leq 6\rho/\sigma^3\sqrt{n}. \]
These inequalities show that for any $\epsilon>0$ there exist $N$
such that $n>N,$

\[\Phi(z_2)-\Phi(z_1)-\epsilon<P(z_1<\sqrt{n}(\frac{\bar{X}_n-\mu}{\sigma}) \leq z_2)< \Phi(z_2)-\Phi(z_1)+\epsilon,\]
for $z_1<z_2 \in \R \cup \{-\infty,\infty\}$.

It is interesting to ask under what conditions $lq_{F_n}$ and
$rq_{F_n}$ tend to $lq_F$ and $rq_F$ as $n \rightarrow \infty$.
Theorem \ref{limit-lq-rq} gives a complete answer to this question.

\begin{theorem} (Quantile Convergence/Divergence Theorem)
\begin{enumerate}[a)]
\item Suppose $rq_F(p)=lq_F(p)$ then
\[rq_{F_n}(p) \rightarrow rq_F(p),\;\;a.s.,\]
and
\[lq_{F_n}(p) \rightarrow lq_F(p),\;\;a.s..\]

\item When $lq_F(p)<rq_F(p)$ then both $rq_{F_n}(p),lq_{F_n}(p)$
diverge almost surely.

\item Suppose $lq_F(p) < rq_F(p)$. Then for every $\epsilon>0$
there exists $N$ such that $n>N,$

\[lq_{F_n}(p),rq_{F_n}(p) \in (lq_F(p)-\epsilon,lq_F(p)] \cup [rq_F(p),rq_F(p)+\epsilon).\]

\item

\[\limsup_{n \rightarrow \infty} lq_{F_n}(p)=\limsup_{n \rightarrow \infty} rq_{F_n}(p)=rq_F(p),\;a.s.,\]

and

\[\liminf_{n \rightarrow \infty} lq_{F_n}(p)=\liminf_{n \rightarrow \infty} rq_{F_n}(p)=rq_F(p),\;a.s..\]

%\item Let $Int(a,b)$ be the interval between $a,b$ (either $(a,b)$ or
%$(b,a)$), then

%\[\lim_{n \rightarrow \infty} P(X \in Int(lq_{F_n}(p),lq_F(p)))=0;a.s.\]

%\[\lim_{n \rightarrow \infty} P(X \in Int(rq_{F_n}(p),rq_F(p)))=0,\;a.s.\]

%In other words:

%\[\lim_{n \rightarrow \infty}\delta_F(lq_{F_n}(p),lq_F(p))=0,\;a.s.\]

%\[\lim_{n \rightarrow \infty}\delta_F(rq_{F_n}(p),rq_F(p))=0,\;a.s.\]
\end{enumerate}
\label{limit-lq-rq}
\end{theorem}

\begin{proof}
\begin{enumerate}[a)]

\item Since, $lq_F(p)=rq_F(p)$, we use $q_F(p)$ to denote both.
Suppose $\epsilon>0$ is given. Then

\[F(q_F(p)-\epsilon)<p  \Rightarrow
F(q_F(p)-\epsilon)=p-\delta_1,\;\delta_1>0,\] and
\[F(q_F(p)+\epsilon)>p \Rightarrow
F(q_F(p)+\epsilon)=p+\delta_2,\;\delta_2>0.\] By the
Glivenko--Cantelli Theorem,

\[F_n(u) \rightarrow F(u)\;\;a.s.,\]
uniformly over $\R$. We conclude that

\[F_n(q_F(p)-\epsilon) \rightarrow
F(q_F(p)-\epsilon)=p-\delta_1,\;\;a.s.,\] and
\[F_n(q_F(p)+\epsilon) \rightarrow
F(q_F(p)+\epsilon)=p+\delta_2,\;\;a.s..\]

Let $\epsilon'=\frac{\min(\delta_1,\delta_2)}{2}$. Pick $N$ such
that for $n>N:$

\[\begin{array}{ll}        p-\delta_1-\epsilon'<F_n(q_F(p)-\epsilon)<
p-\delta_1+\epsilon',\\ p+\delta_2-\epsilon' <F_n(q_F(p)+\epsilon)
<p+\delta_2+\epsilon'. \end{array}
\]
Then
\begin{eqnarray*} F_n(q_F(p)-\epsilon)<
p-\delta_1+\epsilon'<p&\Rightarrow&\\
 lq_{F_n}(p) \geq q_F(p)-\epsilon\; &{\rm and}&\; rq_{F_n}(p) \geq q_F(p)-\epsilon.
 \end{eqnarray*}
 Also
 \begin{eqnarray*}
 p<p+\delta_2-\epsilon' <F_n(q_F(p)+\epsilon)&\Rightarrow&\\
 lq_{F_n}(p)  \leq q_F(p)+\epsilon \; &{\rm and}&\; rq_{F_n}(p)\leq
q_F(p)+\epsilon.
\end{eqnarray*}
 Re-arranging these inequalities we get:

\[q_F(p)-\epsilon \leq lq_{F_n}(p) \leq q_F(p)+\epsilon,\]
and \[q_F(p)-\epsilon \leq rq_{F_n}(p) \leq q_F(p)+\epsilon.\]

\item This needs more development in the sequel and the proof
follows.

\item This also needs more development in the sequel and the proof
follows.

\item If $lq_F(p)=rq_F(p)$ the result follows immediately from
(a). Otherwise suppose $lq_F(p)<rq_F(p)$. Then by (b) $lq_{F_n}(p)$
diverges almost surely. Hence $\limsup lq_{F_n}(p) \neq \liminf
lq_{F_n}(p),\;a.s.\;.$ But by (c), $\forall \epsilon>0,\;\exists
N,\;n>N$

\[lq_{F_n}(p) \in (lq_F(p)-\epsilon,lq_F(p)] \cup [rq_F(p),rq_F(p)+\epsilon).\]
This means that every convergent subsequence of $lq_{F_n}(p)$ has
either limit $lq_F(p)$ or $rq_F(p),\;a.s.$.  Since $\limsup
lq_{F_n}(p) \neq \liminf lq_{F_n}(p),\;a.s.$, we conclude $\limsup
lq_{F_n}(p)=rq_F(p)$ and $\liminf lq_{F_n}(p)=lq_F(p),\;a.s.$.\\
A similar argument works for $rq_{F_n}(p).$
\end{enumerate}
\end{proof}

To investigate the case $lq_F(p)\neq rq_F(p)$ more, we start with
the simplest example namely a fair coin. Suppose $X_1,X_2,\cdots$ an
i.i.d sequence with $P(X_i=-1)=P(X_i=1)=\frac{1}{2}$ and let
$Z_n=\sum_{i=1}^n X_i$. Note that
\[Z_n \leq 0 \Leftrightarrow lq_{F_n}(1/2)=-1,\;\;\;\;\;\;\;\;\;\;\; Z_n > 0 \Leftrightarrow lq_{F_n}(1/2)=1,\]
and
\[Z_n < 0 \Leftrightarrow rq_{F_n}(1/2)=-1,\;\;\;\;\;\;\;\;\;\;\;\;\;\; Z_n \geq 0 \Leftrightarrow rq_{F_n}(1/2)=1.\]
Hence in order to show that $lq_{F_n}(1/2)$ and $lq_{F_n}(1/2)$
diverge almost surely, we only need to show that
$P((Z_n<0\;i.o.)\cap (Z_n>0\;i.o.))=1$.  We start with a theorem
from \cite{breiman}.

\begin{theorem}
Suppose $X_i$ is as above. Then $P(Z_n=0\; i.o.)=1$.
\label{recurring-sum}
\end{theorem}
\begin{proof}
The proof of this theorem in \cite{breiman} uses the Borel-Cantelli
Lemma part 2.
\end{proof}

\begin{theorem}
Suppose, $X_1,X_2,\cdots$ $i.i.d.$ and $P(X_i=-1)=P(X_i=1)=1/2$.
Then $lq_{F_n}(1/2)$ and $rq_{F_n}(1/2)$
 diverge almost surely.
\end{theorem}

\begin{proof}
Suppose, $A=\{Z_n=-1\; i.o.\}$ and $B=\{Z_n=1\; i.o.\}$. It suffices
to show that

\[P(A \cap B )=1.\]
But $\omega \in A \cap B \Rightarrow lq_{F_n}(p)(\omega)=-1,\;
i.o.\;\; \mbox{and}\;\; lq_{F_n}(p)(\omega)=1,\; i.o.$ Hence $
lq_{F_n}(p)(\omega)$ diverges.

Note that $P(A)=P(B)$ by the symmetry of the distribution. Also it
is obvious that both $A$ and $B$ are tail events and so have
probability either zero or one. To prove $P(A\cap B)=1$, it only
suffices to show that $P(A \cup B)>0$. Because then at least one of
$A$ and $B$ has a positive probability, say $A$. \[P(A)>0
\Rightarrow P(A)=1 \Rightarrow P(B)=P(A)=1 \Rightarrow P(A \cap
B)=1.\] Now let $C=\{Z_n=0,\;i.o.\}$. Then $P(C)=1$ by Theorem
\ref{recurring-sum}. If $Z_n(\omega)=0$ then either
$Z_{n+1}(\omega)=1$ or $Z_{n+1}(\omega)=-1$. Hence if
$Z_n(\omega)=0,\; i.o.$ then at least for one of $a=1$ or $a=-1,$
$Z_n(\omega)=a,\; i.o.$. We conclude that $\omega \in A \cup B$.
This shows $C \subset A \cup B \Rightarrow P(A \cup B)=1.$
\end{proof}

To generalize this theorem, suppose $X_1,X_2,\cdots,$ arbitrary
$i.i.d$ process and $lq_F(p)<rq_F(p)$. Define the process

\[Y_i=\begin{cases} 1 & X_i \geq rq_F(p)\\
0 & X_i \leq lq_F(p).\\\end{cases}\] (Note that
$P(lq_X(p)<X<rq_X(p))=0$.) Then the sequence $Y_1,Y_2,\cdots$ is
$i.i.d.$, $P(Y_i=0)=p$ and $P(Y_i=1)=1-p$. Also note that

\[lq_{F_{n,Y}}(p) \; \mbox{diverges a.s. } \Rightarrow lq_{F_{n,X}}(p)\; \mbox{diverges a.s.}\]
Hence to prove the theorem in general it suffices to prove the
theorem  for the $Y_i$ process. However, we first prove a lemma that
we need in the proof.

\begin{lemma}
Let $Y_1,Y_2,\cdots$ $i.i.d$ with $P(Y_i=0)=p=1-q>0$ and
$P(Y_i=1)=1-p=q>0$. Let $S_n=\sum_{i=1}^n
Y_i,\;0<\alpha,\;\;k\in\N$. Then there exists a transformation
$\phi(k)$ (to $\N$) such that
\[P(S_{\phi(k)}-\phi(k)q<-k)>1/2-\alpha,\]
\[P(S_{\phi(k)}-\phi(k)q>k)>1/2-\alpha.\]
\label{lemma-sum-deviations}
\end{lemma}

\noindent{\bf Remark.} For $\alpha=1/4$, we get
\[P(S_{\phi(k)}-\phi(k)q<-k)>1/4,\]
\[P(S_{\phi(k)}-\phi(k)q>k)>1/4.\]

\begin{proof}
Since the first three moments of $Y_i$ are finite
($E(Y_i)=q,E(|Y_i-q|^2)=q(1-q)=\sigma,E(|Y_i-q|^3)=q^3(1-q)+(1-q)^3q=\rho$),
we can apply the Berry-Esseen theorem to
$\sqrt{n}\frac{\bar{Y}_n-\mu}{\sigma}$. By a corollary of that
theorem, for $\frac{\alpha}{2}>0$ there exists an $N_1$ such that

\[1-\Phi(z)-\frac{\alpha}{2}<P(\sqrt{n}\frac{\bar{Y}_n-\mu}{\sigma}>z)<1-\Phi(z)+\frac{\alpha}{2},\]
and
\[\Phi(z)-\frac{\alpha}{2}<P(\sqrt{n}\frac{\bar{Y}_n-\mu}{\sigma}<-z)<\Phi(z)+\frac{\alpha}{2},\]
for all $z \in \R$ and $n>N_1$. Now for the given integer $k$ pick
$N_2$ such that

\[\frac{1}{2}-\frac{\alpha}{2}<\Phi(\frac{k}{\sigma \sqrt{N_2}})<\frac{1}{2}+\frac{\alpha}{2}.\]
This is possible because $\Phi$ is continuous and $\Phi(0)=1/2$. Now
let \[\phi(k)=\max\{N_1,N_2\},\;z=\frac{k}{\sigma \sqrt{\phi(k)}}.\]
Then since $\phi(k) \geq N_1$

\[P(\sqrt{\phi(k)}\frac{\bar{Y}_{\phi(k)}-\mu}{\sigma}>z)>1-\Phi(z)-\frac{\alpha}{2}>1/2-\alpha,\]
and
\[P(\sqrt{\phi(k)}\frac{\bar{Y}_{\phi(k)}-\mu}{\sigma}<-z)>\Phi(z)-\frac{\alpha}{2}>1/2-\alpha.\]
These two inequalities are equivalent to
\[P((S_{\phi(k)}-\phi(k)q)<-k)>1/2-\alpha,\]
and
\[P((S_{\phi(k)}-\phi(k)q)>k)>1/2-\alpha.\]

If we put $\alpha=1/4$, we get
\[P((S_{\phi(k)}-\phi(k)q)<-k)>1/4,\]
and
\[P((S_{\phi(k)}-\phi(k))q>k)>1/4.\]

\end{proof}

We are now ready to prove Part b) of Theorem \ref{limit-lq-rq}.

\begin{proof} [Theorem \ref{limit-lq-rq}, Part b)]

For the process $\{Y_i\}$ as defined above, let $n_1=1, m_k=n_k +
\phi(n_k)$ and $n_{k+1}=m_k+\phi(m_k)$. Then define

\[D_k=(Y_{n_k+1}+\cdots+Y_{m_k}-(m_k-n_k)q<-n_k),\]
\[E_k=(Y_{m_k+1}+\cdots+Y_{n_{k+1}}-(n_{k+1}-m_k)q>m_k),\]
\[C_K=D_k \cap E_k.\]
Since $\{C_k\}$ involve non-overlapping subsequences of $Y_s$, they
are independent events. Also $D_k$ and $E_k$ are independent. Now
note that

\begin{align*}
Y_{n_k+1}+\cdots+Y_{m_k}-(m_k-n_k)q<-n_k \Rightarrow\\
Y_1+\cdots+Y_{m_k}<-n_k+(m_k-n_k)q+n_k  \Rightarrow\\
 \bar{Y}_{m_k}<\frac{m_k-n_k}{m_k}q<q \Rightarrow \\
lq_{F_{n,Y}}(p)=rq_{F_{n,Y}}=0 \Rightarrow\\
 \{C_k,\;\;i.o.\}\subset
\{lq_{F_{n,Y}}(p)=rq_{F_{n,Y}}=0,\;\;i.o.\}.
\end{align*}
Similarly,

\begin{align*}
Y_{m_k+1}+\cdots+Y_{n_{k+1}}-(n_{k+1}-m_k)q>m_k\\ \Rightarrow
Y_1+\cdots+Y_{n_{k+1}}>(n_{k+1}-m_k)q+m_k\\
  \Rightarrow \bar{Y}_{n_{k+1}}>\frac{m_k+(n_{k+1}-m_k)q}{n_{k+1}}>q=1-p\\ \Rightarrow
lq_{F_{n,Y}}(p)=rq_{F_{n,Y}}(p)=1\\\Rightarrow
\{C_k,\;\;i.o.\}\subset
\{lq_{F_{n,Y}}(p)=rq_{F_{n,Y}}(p)=1,\;\;i.o.\}.
\end{align*}
Let us compute the probability of $C_k$:

\begin{align*}
P(C_k)=\;\;\;\;\;\;\;\;\;\;\;\;\;\;\;\;\;\;\;\;\;\;\;\;\;\;\;\;\;\;\;\;\;\;\;\;\;\\
P(Y_{n_k+1}+\cdots+Y_{m_k}-(m_k-n_k)q<-n_k) \times\\
P(Y_{m_k+1}+\cdots+Y_{n_{k+1}}-(n_{k+1}-m_k)q>m_k)=\\
P(Y_{1}+\cdots+Y_{\phi(n_k)}-\phi(n_k)q<-n_k) \times
\\P(Y_{1}+\cdots+Y_{\phi(m_k)}-\phi(m_k)q>m_k)>1/4.1/4=1/16.
\end{align*}
We conclude that
\[\sum_{k=1}^{\infty}P(C_k)=\infty.\]
By the Borel-Cantelli Lemma, $P(C_k,\;\;i.o.)=1$. We conclude that
\[P(lq_{F_{n,Y}}(p)=rq_{F_{n,Y}}(p)=0,\;\;i.o.)=1,\] and \[P(lq_{F_{n,Y}}(p)=rq_{F_{n,Y}}(p)=1,\;\;i.o.)=1.\]
Hence,
\[P(\{lq_{F_{n,Y}}(p)=rq_{F_{n,Y}}(p)=0,\;\;i.o.\} \cap \{lq_{F_{n,Y}}(p)=rq_{F_{n,Y}}(p)=1,\;\;i.o.\})=1.\]
\end{proof}

\begin{figure}
\centering
\includegraphics[width=0.7\textwidth] {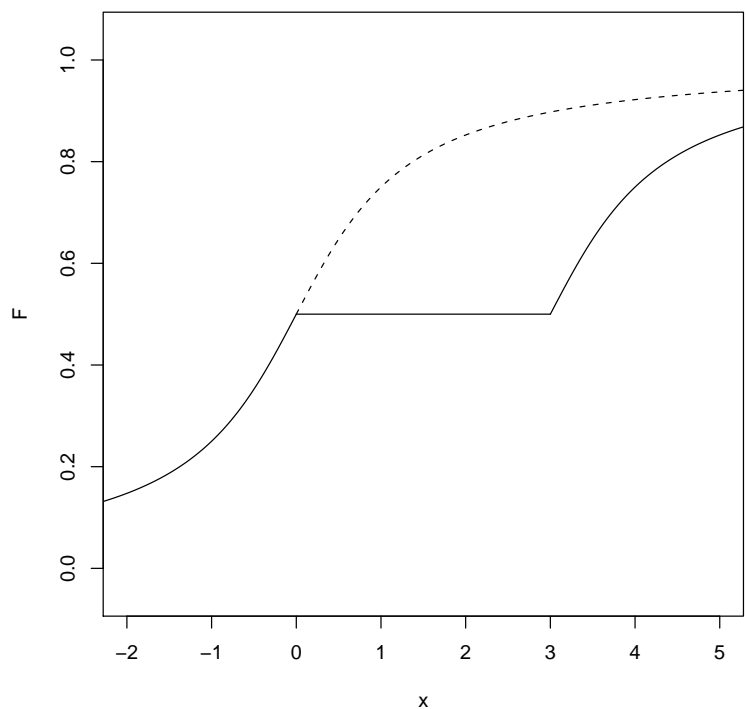}
 \caption{The solid line is the distribution function of $\{X_i\}$. Note that for the distribution of the $X_i$
  and $p=0.5$, $lq_{F_X}(p)=0,rq_{F_X}(p)=3$.
Let $h=rq(p)-lq(p)=3.$ The dotted line is the distribution function
of the $\{Y_i\}$ which coincides with that of $\{X_i\}$ to the left
of $lq_{F_X}(p)$ and is a backward shift of 3 units for values
greater than $rq_{F_X}(p)$. Note that for the $\{Y_i\}$,
$lq_{F_Y}(p)=rq_{F_Y}(p)=1.$}
 \label{proof-dis-fcn.ps}
\end{figure}

\begin{proof} (Theorem \ref{limit-lq-rq}, part (c))\\
Suppose that $rq_F(p)=x_1 \neq lq_F(p)=x_2$ and $a$ is an arbitrary
real number. Let $h=x_2-x_1$. We define a new chain $Y$ as follows:

\[Y_i=\begin{cases}X_i  & X_i \leq lq_{F_X}(p)\\
X_i-h & X_i \geq rq_{F_X}(p). \end{cases}\] (See Figure
\ref{proof-dis-fcn.ps}.) Then $Y_1,Y_2,\cdots$ is an $i.i.d$ sample.
We drop the index $i$ from $Y_i$ and $X_i$ in the following for
simplicity and since the $Y_i$ (as well as the $X_i$) are
identically distributed. We claim

\[lq_{F_Y}Y(p)=rq_{F_Y}(p)=lq_{F_X}(p).\]
To prove $lq_{F_Y}(p)=lq_{F_X}(p)$, note that
\[F_Y(lq_{F_X}(p))=P(Y \leq lq_{F_X}(p)) \geq P(X \leq lq_{F_X}(p)) \geq p \Rightarrow lq_{F_Y}(p) \leq lq_{F_X}(p).\]
(The first inequality is because $Y \leq X$.) Moreover for any $y <
lq_{F_X}(p)$, $F_Y(y)=F_X(y)<p$. (Since  $X,Y<lq_{F_X}(p)
\Rightarrow X=Y$.) Hence $lq_{F_Y}(p) \geq lq_{F_X}(p)$ and we are
done. To show $rq_{F_Y}(p)=lq_{F_X}(p)$, note that $rq_{F_Y}(p)\geq
lq_{F_Y}(p)=lq_{F_X}(p)$. It only remains to show that $rq_{F_Y}(p)
\leq lq_{F_X}(p)$. Suppose $y>lq_{F_X}(p)$ and let
$\delta=y-lq_{F_X}(p)>0$. First note that \begin{eqnarray*} P(\{Y
\leq lq_{F_X}(p)+\delta\})&=&\\
P(\{Y \leq lq_{F_X}(p)+\delta \mbox{ and } X \geq rq_{F_X}(p)\}
&\cup&\\\{Y \leq lq_{F_X}(p)+\delta \mbox{ and } X \leq
lq_{F_X}(p)\})&=&\\
P(\{X-h \leq lq_{F_X}(p)+\delta \mbox{ and } X \geq rq_{F_X}(p)\}
&\cup&\\\{X \leq lq_{F_X}(p)+\delta \mbox{ and } X \leq
lq_{F_X}(p)\})&=&\\
P(\{rq_{F_X}(p) \leq X \leq rq_{F_X}(p)+\delta\} \cup \{X \leq
lq_{F_X}(p)\})&=&\\
P(\{X \leq rq_{F_X}(p)+\delta\}).&&\end{eqnarray*} Hence,

\[F_Y(y)=P(Y \leq lq_{F_X}(p)+\delta)=P(X \leq rq_{F_X}(p)+\delta)>p \Rightarrow\] \[rq_{F_Y}(p) \leq y,\; \forall y>lq_{F_X}(p).\]
We conclude that $rq_{F_Y}(p)\leq lq_{F_Y}(p).$\\
To complete the proof of part (c) observe that for every
$\epsilon>0$, we may suppose that $lq_{F_{n,Y}}(p) \in
(q_{F_Y}(p)-\epsilon,q_{F_Y}(p)+\epsilon).$ Then

\begin{equation}
lq_{F_n,X}(p),rq_{F_{n,X}}(p) \in
(lq_{F_X}(p)-\epsilon,rq_{F_X}(p)+\epsilon).
\label{claim-quatile-limit-part-c}
\end{equation}
This is because from $lq_{F_{n,Y}}(p)
\in(q_{F_Y}(p)-\epsilon,q_{F_Y}(p)+\epsilon)$, we may conclude that

\[F_{n,Y}(q_{F_Y}(p)+\epsilon)>p \Rightarrow F_{n,X}(rq_{F_X}(p)+\epsilon)>p\Rightarrow\] \[lq_{F_{n,X}}(p),rq_{F_{n,X}}(p) <rq_{F_X}(p)+\epsilon, \]
and
\[F_{n,Y}(q_{F_Y}(p)-\epsilon)<p \Rightarrow F_{n_X}(lq_{F_X}(p)-\epsilon)<p\Rightarrow\] \[lq_{F_{n,X}}(p),rq_{F_{n,X}}(p) >lq_{F_X}(p)-\epsilon.\]
But by part (a) of Theorem \ref{limit-lq-rq}, $lq_{F_{n,Y}}(p)
\rightarrow q_{F_Y}(p)$ and $rq_{F_{n,Y}}(p) \rightarrow
q_{F_Y}(p)$. Hence for given $\epsilon>0$ there exists an integer
$N$ such that for any $n>N,$ $lq_{F_{n,Y}}(p) \in
(q_{F_Y}(p)-\epsilon,q_{F,Y}(p)+\epsilon)$. By
(\ref{claim-quatile-limit-part-c}), we have shown that for every
$\epsilon>0$ there exists $N$ such that for every $n>N$

\[q_{F_{n,X}}(p),rq_{F_{n,X}}(p) \in (lq_{F_X}(p)-\epsilon,rq_{F_X}(p)+\epsilon),\]
since
\[P(X_i \in (lq_{F_X}(p),rq_{F_X}(p))\;\mbox{for some}\; i \in
\N)=0.\] We can conclude that

\[P(lq_{F_{n,X}}(p)\in
(lq_{F_X}(p),rq_{F_X}(p))\;\mbox{for some}\; i \in \N)=0\] and
\[P(rq_{F_{n,X}}(p)\in (lq_{F_X}(p),rq_{F_X}(p))\;\mbox{for some}\; i \in
\N)=0.\] Hence with probability 1

\[q_{F_{n,X}}(p),rq_{F_{n,X}}(p) \in (lq_{F_X}(p)-\epsilon,lq_{F_X}(p)] \cup [rq_{F_X}(p),rq_{F_X}(p)+\epsilon).\]
\end{proof}

\newpage

\bibliographystyle{plain}
\bibliography{mybibreza}
\end{document}